\documentclass[11pt]{article}
\usepackage{amssymb,amsfonts,amsthm,ifpdf}
\usepackage{color}
\newtheorem{thm}{Theorem}[section]
\newtheorem{cor}[thm]{Corollary}

\newtheorem{prop}[thm]{Proposition}

\def\eop{\hfill\rule{2.5mm}{2.5mm}}
\def\pf{\par\smallbreak\noindent {\bf Proof.} \ }

\begin{document}

\title{
{\textbf{\Large{A close look at the entropy numbers of the unit ball of the Reproducing Hilbert Space of isotropic positive definite kernels}}} \vspace{-4pt}
\author{
\textsc{T. Jord\~{a}o}\,\footnote{FAPESP, grant no 2022/11032-0.}
\,\,\&\,\,
\textsc{K. Gonzalez} \footnote{Partially supported by CAPES, finance code 001.}}
}
\date{}

\maketitle 

\vspace{-30pt}
\begin{center}
\parbox{13 cm}{{\small {\bf Abstract:} We present upper and lower bounds for the covering numbers, with accurate explicit constants, of the unit ball for two general classes of Reproducing Kernel Hilbert Space (RKHS) on the unit sphere of $\mathbb{R}^{d+1}$.\ In both classes, the RKHS is generated by an isotropic continuous positive definite kernel and the bounds we present carry precise information about the asymptotic constants, depending on the dimension of the sphere and the monotonic behavior of the Schoenberg/Fourier coefficients of the isotropic kernel. }}
\end{center}

\thispagestyle{empty}

\section{Introduction}
\label{intro}

The $\epsilon$-entropy introduced by Kolmogorov (\cite{kolmogorov}), or the covering numbers, on manifolds plays an important role in important areas such as kernel-based learning algorithms, Gaussian process, and density estimation (\cite{devroyebook,Li,Quang,WSS}).\ The asymptotic behaviour or bounds for the covering numbers turned a very important tool in order to estimate the probabilistic error of the statistical nature of the observations from which the algorithms employed are learning from (\cite{smale,zhou1}).\ The explicit constants and the sharpness included in the process of finding the bounds for the covering numbers are also very important as recently illustrated in \cite{kuhn22, ingo}.\ The author in \cite{ingo} consider the Gaussian kernel and recovers a previous well-know asymptotic behaviour of the covering number obtained in \cite{Kuhn}, but preserving the constants expressed in terms of the kernel and the dimension of its domain. 

In this paper we provide precise estimates for the covering numbers of the unit ball of two general classes of Reproducing kernel Hilbert Space (RKHS) on the unit sphere.\ We follow the path designed in \cite{Gonzalez}, taking in account the generality of the context and the characterization of isotropic (radially symmetric) continuous positive definite kernels.\ We give explicit constants and sharpness for the bounds for the covering numbers, as in \cite{ingo} for the Gaussian Kernel.\ Isotropic positive definite kernels on the sphere are fundamental in spacial statistics (see \cite{PorcuBerg,Montserrat}) and the Schoenberg representation (\cite{Gneiting,Schoenberg}) is widely considered on those scenarios.\ This representation has direct relation with the Fourier series expansion of the kernels and we will explore it fundamentally.\ The relation exists more generally even for non radially symmetric kernels as can be seen in \cite{JFAA22}.\ The Gaussian Kernel is the classical and useful example of isotropic positive kernel (see, for example, \cite{Quang}), other examples such as the Powered exponential, Matérn, and Generalized Cauchy kernels can be found in \cite{Gneiting}.

In order to present our results we first present the basic and general setting.\ Let $\mathbb{S}^d$ be the unit sphere centered at the origin of the $(d+1)$-dimensional euclidean space, for $d\geq 2$.\ We consider $\mathbb{S}^d$ endowed with the induced Lebesgue measure (surface measure) denoted by $\sigma_d$.\ We work with the standard Hilbert space $(L^2(\mathbb{S}^d), \langle\, \cdot\, , \, \cdot \, \rangle_2)$ of square integrable functions on $\mathbb{S}^d$, $L^2(\mathbb{S}^d):=L^2(\mathbb{S}^d, \sigma_d)$, and with the induced norm given by 
$$
\| f \|_2 = \left(\frac{1}{\sigma_d(\mathbb{S}^d)}\int_{\mathbb{S}^d}|f(x)|^2d\sigma_d(x)\right)^{1/2}, \quad f\in L^2(\mathbb{S}^d).
$$
For $k=0,1,\ldots$, we write $H^k (\mathbb{S}^d)$ the space of spherical harmonics of degree $k$ in $d+1$ variables and $\tau_d^k:= \mbox{dim}\,H^k(\mathbb{S}^d)$.\ If $k\neq l$, then $H^k (\mathbb{S}^d )\perp H^l (\mathbb{S}^d)$ and we write $\{S_{k,1}, S_{k,2} ,..., S_{k,\tau_d^k}\}$ for the orthonormal basis of $H^k(\mathbb{S}^d)$.\ For each $f\in L^2(\mathbb{S}^d)$ there exists a unique \emph{Fourier series expansion} of $f$, given by 
$$
f\sim\sum_{k=0}^{\infty}\sum_{j=1}^{\tau_d^k}\langle f, S_{k,j}\rangle_2 S_{k,j}.
$$
The relation between Fourier series expansions of spherical functions and the so-called Schoenberg series representation is given in terms of the well-known orthogonal polynomials, the Legendre polynomials.\ If we write $P^d_{k}$ for the Legendre polynomial of degree $k$ and $d+1$ variables, then we have the following useful identity 
\begin{equation}\label{addition}
\sum_{j=1}^{\tau_d^k}S_{k,j}(x)S_{k,j}(y)=\tau_d^kP^{d}_{k}(x\cdot y), \quad x,y\in \mathbb{S}^d, \quad k=0,1,\ldots,
\end{equation}
named the Addition Formula (see \cite[Section 2.3]{morimoto}).\ This identity promptly implies $P^{d}_{k}(1)=1$, and 
$$
|P^{d}_{k}(t)|\leq 1, \quad \mbox{for} \,\, -1\leq t\leq 1, \quad k=0,1,\ldots.
$$

In this paper we will work with an isotropic or, equivalently, a zonal kernel. A function $K:\mathbb{S}^d\times \mathbb{S}^d\to \mathbb{R}$ is an {\emph{isotropic kernel} if $K(x,y)=f(x\cdot y)$, for all $x,y\in  \mathbb{S}^d$ and for some continuous function $f:[-1,1]\longrightarrow \mathbb{R}$.\ We will deal with $K$ a \emph{positive definite} kernel,  it means that for any $n\in\mathbb{N}$, $x_1,\ldots,x_n \in \mathbb{S}^d$ and $c_1, \ldots,c_n\in \mathbb{R}$ it holds
$$
\sum_{i=1}^{n}\sum_{j=1}^{n}c_{i}c_{j}K(x_{i},x_{j})\geq 0,
$$
if the inequality above is strict, then we say that the kernel $K$ is {\emph{strictly positive definite}}.\ The classical Schoenberg's theory (\cite[Theorem 1.1]{peron} and \cite{PorcuBerg}) asserts that an isotropic positive definite kernel $K$ can be written as a series expansion in terms of the Legendre polynomial, with summable sequence of non-negative real numbers $\{a_k\}$ as Schoenberg coefficients, it means
\begin{equation}\label{WKlegendre0}
K(x,y)=\sum_{k=0}^\infty a_k P_k^d (x\cdot y), \quad x,y\in \mathbb{S}^d.
\end{equation}
The Addition formula (\ref{addition}), based on the Schoenberg representation of $K$, promptly gives us a Fourier series representation for $K$ as follows
\begin{equation}\label{WKlegendre}
K(x,y)=\sum_{k=0}^\infty a_k P_k^d (x\cdot y) = \sum_{k=0}^{\infty} \frac{a_k}{\tau_d^k} \sum_{j=1}^{\tau_d^k} S_{k,j} (x) S_{k,j} (y), \quad x,y\in \mathbb{S}^d.
\end{equation}
We assume that $a_k\downarrow 0$, as $k\to\infty$.
If $K$ is strictly positive definite, it is well known that $\{n\in\mathbb{N} : a_n>0\}$ is infinity and contain infinitely many odd integers as well as many even integers (\cite{PD2}).\ 

The general Aronszajn theory about RKHSs (\cite{aron}) ensures that for a positive definite kernel $K$, there exists a unique Hilbert space $(\mathcal{H}_{K},\langle \,\cdot\, , \cdot \, \rangle_K)$ of functions on $\mathbb{S}^d$ satisfying:\\
i. $K(x,\, \cdot\, )\in \mathcal{H}_{K}$ for all $x\in \mathbb{S}^d$;\\
ii. span$\{K(\, \cdot\, ,x) : x\in X\}$  is dense in $\mathcal{H}_{K}$;\\
iii. (Reproducing property) $f(x)=\langle{f},{K(\cdot,x)}\rangle_{K}$, for all $x\in\mathbb{S}^d$ and $f\in \mathcal{H}_{K}$.\\
From the reproducing property, since $K$ is a continuous kernel, we have that $\mathcal{H}_{K}\hookrightarrow C(\mathbb{S}^d)$, with $C(\mathbb{S}^d)$ the space of continuous functions with the supremum norm $\|\,\cdot \,\|_{\infty}$.\ Precisely, this embedding $I_{K}:\mathcal{H}_K \longrightarrow C(\mathbb{S}^d)$ is to one we are interested on the estimates for the covering numbers.\ 

The covering numbers of the operator $I_K$ are defined in terms of the covering numbers of subsets of metric spaces.\ If $A$ is a subset of a metric space $X$ and $\epsilon>0$, then the covering number $\mathcal{C} (\epsilon, A)$ is the minimal number of balls of $X$ with radius $\epsilon$ which covers $A$.\ We denote by $B_K$ the unit ball in $\mathcal{H}_{K}$ and the \emph{covering numbers} of the embedding $I_K$ are defined as follows
$$
\mathcal{C} (\epsilon, I_K):= \mathcal{C}(\epsilon,I_K(B_K)), \quad \epsilon>0.
$$

The first class of kernels we consider contains just infinitely many times differentiable kernels.\ This conclusion is an easy application of the criteria in \cite[Theorem 2.45]{morimoto} for kernels under the hypotheses of Theorem \ref{coveringestimatesgeo}.\ This class is given by those kernels having the sequence of Schoenberg coefficients decreasing no more slowly than a geometric progression (\cite{bary}) and the Gaussian kernel is an example of such kernel.\ For the Gaussian kernel the asymptotic analysis for the covering number can be found in \cite{Kuhn}.\ Our result reads as follows.\ 

\begin{thm}\label{coveringestimatesgeo} Let $K$ be an isotropic positive definite kernel given by the series expansion (\ref{WKlegendre}).\ If there exist $\delta, \theta\in (0,1)$ such that  $0<\delta\,a_{k-1}<a_k \leq \theta\, a_{k-1}$, for $k=1,2\ldots$, then
$$
\frac{ (2d)^d }{d! (d + \ln (1/\delta)-1)^d (d+1)^{d+1}} \leq \liminf_{\epsilon\to 0^+}\frac{\ln(\mathcal{C}(\epsilon, I_K ))}{[\ln \left(1/\epsilon \right)]^{d+1}} \leq \limsup_{\epsilon\to 0^+}\frac{\ln(\mathcal{C}(\epsilon, I_K ))}{[\ln \left(1/\epsilon \right)]^{d+1}} \leq \frac{ 2^{d+2} }{ d! \left({\ln\left( {1}/{\theta}\right)}\right)^d} .
$$
\end{thm}

In the second result, we consider kernels having the sequence of Schoenberg/Fouriers coefficients with asymptotic behavior as a harmonic sequence.\ The asymptotic analysis for the covering numbers of these kernels were considered in \cite{Azevedo} on a general manifold, including the spherical setting, but with no concern about the constants in the bounds.\ Also, the technique employed in order to obtain the lower bounds is slightly different from here.

\begin{thm}\label{coveringestimateshar} Let $K$ an isotropic positive definite kernel given by the series expansion (\ref{WKlegendre}).\ If there exist $\beta>0$, $\rho>1$, and $c_1,c_2>0$ such that $c_2 k^{-\rho}\leq a_k\leq c_1 k^{-\beta-d}$, then
\begin{itemize}
    \item[(I)]     $$
\limsup_{\epsilon\to 0^+}\frac{\ln (\mathcal{C} (\epsilon, I_K))}{(1/\epsilon)^{2d/(\beta+d-1)}\ln (1/\epsilon)}\leq \frac{4}{d!} \left(\frac{4 c_1 }{\beta+d-1}\right)^{d/(\beta+d-1)},
$$
    \item[(II)] $$ \ln{\sqrt{a_0}}+ \frac{2(\rho+d-1)}{e\, d^2 (d-1)!}  \left(\frac{c_2}{2 } \right)^{d/(\rho+d-1)} \leq \liminf_{\epsilon\to 0}\frac{\ln(\mathcal{C}(\epsilon, I_K ))}{(1/\epsilon)^{d/2(d+\rho-1)}}. 
$$
    \end{itemize}
\end{thm}

In the rest of the paper let two given functions $f,g: (0,\infty) \longrightarrow \mathbb{R_+}$, we write $f(\epsilon) \asymp g(\epsilon)$ for the weak equivalence and it stands for $f(\epsilon)=O(g(\epsilon))$ and $g(\epsilon)=O(f(\epsilon))$, as $\epsilon\to 0$.\ For the strong equivalence we write $f(\epsilon) \approx g(\epsilon)$ and it means that $f(\epsilon)/g(\epsilon) \to 1$, as $\epsilon\to 0$.\ 

The paper is organized as follows.\ In Section \ref{RKHS}, we present some general facts about RKHSs and the characterization of them for isotropic positive definite kernels $K$ represented by the series expansion (\ref{WKlegendre}).\ This gives us an orthonormal basis for $\mathcal{H}_{K}$ and we explore the operator norm of the embedding $\mathcal{H}_K \hookrightarrow C(\mathbb{S}^d)$, and related ones.\ This will be fundamentally used in Section \ref{estimates} where the upper and the lower bounds for the covering numbers are obtained in our main propositions.\ The proofs of Theorem \ref{coveringestimatesgeo} and Theorem \ref{coveringestimateshar} are presented in Subsection \ref{estimates1} and Subsection \ref{estimates2}, respectively.\ In Section \ref{consequences}, we explore the case $0<a_k\leq \theta^k$, for $0<\theta<1$ and, consequences of our main results and final comments are presented.\ We finish the paper with Section 5, with some comments about the assumption made in Theorem \ref{thmlowerhar} and the basic decay of the coefficients of isotropic positive definite kernels, and we close the paper with the acknowledgements.   

\section{The RKHS and the embedding $\mathcal{H}_K \hookrightarrow C(\mathbb{S}^d)$}\label{RKHS}

The proofs of the theorems apply general estimates for the covering numbers given em terms of operator norm and for operators with finite dimensional range.\ Then, we present an ortonormal basis for $\mathcal{H}_K$, the operator norm of the embedding $\mathcal{H}_K \hookrightarrow C(\mathbb{S}^d)$, and the operator norm of the embeddings of projections onto finite dimensional subspaces of $\mathcal{H}_K$. 

The following proposition presents the characterization of the RKHS in terms of Fourier series expansions.

\begin{prop}\label{RKHS} If $K$ is an isotropic positive definite kernel represented by the series expansion (\ref{WKlegendre}), then
$$
\mathcal{H}_K = \left\{g:\mathbb{S}^d \longrightarrow \mathbb{R} :
 g(x)=\sum_{k=0}^{\infty} \sqrt{\frac{a_k }{\tau_d^k}} \sum_{j=1}^{\tau_d^k} c_{k}^{j} S_{k,j} (x), \quad x\in\mathbb{S}^d, \, \{c_{k}^{j}\} \in \ell^2 , \, j=1,..,\tau_d^k \right\}
$$
endowed the inner product
$$
\langle  g,h \rangle_{K} =\sum_{k=0}^{\infty} \sum_{j=1}^{\tau_d^k} c_{k}^{j} d_{k}^{j}, \quad g,h \in \mathcal{H}_K,
$$
where
$$
g(x)= \sum_{k=0}^{\infty} \sqrt{\frac{a_k}{\tau_d^k}} \sum_{j=1}^{\tau_d^k} c_{k}^{j} S_{k,j} (x) \quad \mbox{and}\quad h(x)= \sum_{k=0}^{\infty} \sqrt{\frac{a_k}{\tau_d^k}} \sum_{j=1}^{\tau_d^k} d_{k}^{j} S_{k,j} (x), \quad x\in \mathbb{S}^d.
$$ 
\end{prop}

Due the unicity of the RKHS, this fact can be easily proved just by showing that properties from i) to iii) in Section \ref{intro} are satisfied for $\mathcal{H}_K$.\ This is standard in theory of RKHSs and can be seen in \cite[Theorem 2.2]{Gonzalez}, for example.

We consider $P_m$ the orthogonal projections onto the subspace $\mathcal{V}_m$ of $\mathcal{H}_K$ with the basis
$$
\left\{\sqrt{\frac{a_k}{\tau_d^k}}S_{k,j}: \quad j=1,..,\tau_d^k, \,\, k=0,1,\ldots, m\right\},
$$
and $P_m^s$ the orthogonal projection onto the subspace and $\mathcal{V}_m^s$, the orthogonal complement of $\mathcal{V}_m$, with orthonormal complete system given by 
$$
\left\{\sqrt{\frac{a_k}{\tau_d^k}}S_{k,j}: \quad j=1,..,\tau_d^k, \,\, k=m+1, m+2, \ldots\right\}. 
$$
Then, we can write
$$
f=P_m (f) + P_m^s(f), \quad f\in\mathcal{H}_{K}= \bigoplus_{k=0}^{\infty}\mathcal{V}_k, \quad m=1,2,\ldots,
$$
and we consider the natural embeddings $\mathcal{V}_m \hookrightarrow C(\mathbb{S}^d)$ and $\mathcal{V}_m^s \hookrightarrow C(\mathbb{S}^d)$, keeping the notation $P_m$ and $P_m^s$ for them, respectively.\ From the well known (\cite[Section 2.1]{morimoto}) facts
\begin{equation}\label{dim}
\tau_d^k = \sum_{j=0}^{k} \tau_{d-1}^j=\frac{(2k+d-1)(k+d-2)!}{k!(d-1)!}, \quad k=0,1,\ldots,
\end{equation}
we observe that 
$$
\dim \mathcal{V}_m \leq \sum_{k=0}^m \tau_d^k = \tau_{d+1}^m = \frac{(2m+d)(m+d-1)!}{m!d!}, \quad m=1,2,\ldots.
$$

Immediate consequence of the previous considerations is the following.

\begin{cor}\label{dimVm} Let $K$ be an isotropic positive definite kernel represented by the series expansion (\ref{WKlegendre}). If $a_k>0$, for $k=0,1,2,\ldots$, then
$$
\dim \mathcal{V}_m = \tau_{d+1}^m, \quad m=1,2,\ldots,
$$
and there exists $m_0\in \mathbb{N}$ such that
\begin{equation}\label{dimapproxm}
\frac{1}{d!} <\frac{\dim \mathcal{V}_m}{m^{d}}< \frac{4}{d!}, \quad  m=m_0,m_0+1,\ldots.
\end{equation} 
\end{cor}

An application of Stirling's formula for factorial in formula (\ref{dim}), imply the following useful strong equivalence
\begin{equation}\label{dimapprox}
\frac{\tau_d^k}{k^{d-1}}\approx \frac{2}{(d-1)!}, \quad \mbox{as \,\, $k\to\infty$},
\end{equation}
and this is enough to prove inequality (\ref{dimapproxm}).

The operator norm of the previous embeddings in terms of the Fourier coefficients of the kernel are given in the following result.

\begin{prop}\label{thmNorms} Let $K$ an isotropic positive definite kernel represented by the series expansion (\ref{WKlegendre}).\ For $m=1,2,\ldots$, the following embeddings
$$
I_K : \mathcal{H}_K \longrightarrow C(\mathbb{S}^d), \quad P_m: \mathcal{V}_m \longrightarrow C(\mathbb{S}^d), \quad \mbox{and} \quad P_m^s: \mathcal{V}_m^s \longrightarrow C(\mathbb{S}^d),
$$
satisfy respectively,
$$
\| I_K \|^2 = \sum_{k=0}^{\infty}a_k, \quad \|P_m\|^2 = \sum_{k=0}^{m}a_k, \quad \mbox{and} \quad
\|P_m^s\|^2 = \sum_{k=m+1}^{\infty}a_k.
$$
\end{prop}


\pf We will show the first part of the theorem just for $P_m$, for $m=1,\ldots$, since the others cases are completely analogous.\ We first observe that
$$
\Vert P_m \Vert^2 = \sup_{g\in\mathcal{V}_m, \| g \|_K =1} \sup_{x\in \mathbb{S}^d} \| \langle   g, K(x,\, \cdot \,) \rangle_K  \| ^2,
$$
and for $g\in\mathcal{V}_m$, it holds
$$
\langle  g, K_M(x,\, \cdot \,) \rangle_K = \langle  g, P_m(K_M(x,\, \cdot \,)) \rangle_K, \quad x\in\mathbb{S}^d. 
$$
The reproducing property of $K$ and Cauchy-Schwarz inequality, imply  
\begin{eqnarray*}
\Vert P_m \Vert^2 & \leq & \sup_{\| g \|_K =1, g\in\mathcal{V}_m} \sup_{x\in \mathbb{S}^d} \| g\|_K ^2 \| K_M (x,\cdot)\|_K^2 \\
   & = & P_m(K_M (x,\cdot))(x).
\end{eqnarray*}
If we consider $x\in\mathbb{S}^d$, then
$$
g_x(w) := \frac{P_m(K(x, \,\cdot\,))(w)}{\|P_m(K(x, \,\cdot\,))\|_K}, \quad w\in\mathbb{S}^d,
$$
is such that $g_x\in \mathcal{V}_m$ and $\|g_x\|_K=1$. Additionally, it is not hard to see that
$$
\| P_m \|^2\geq \|P_m(g_x)\|_{\infty}^2 = \sup_{w\in\mathbb{S}^d} \left|\frac{P_m(K(x, \,\cdot\,))(w)}{\|P_m(K(x, \,\cdot\,))\|_K}\right|^2 = [g_x(x)]^2,
$$
and
$$
[g_x(x)]^2 = \frac{[P_m(K(x, \,\cdot\,))(x)]^2}{\|P_m(K(x, \,\cdot\,))\|^2_K} = P_m(K(x, \,\cdot\,))(x),
$$
Then, \[ \Vert P_m \Vert^2 =\sup_{x\in\mathbb{S}^d}|P_m(K(x,\,\cdot\,))(x)|, \quad m=0,1,\ldots.\]
The definition of $K$ given in (\ref{WKlegendre}) and an application of the Addition formula (\ref{addition}) finishes the proof.\ \eop

\section{Estimates for the covering number}\label{estimates}

We present the basic and general properties of the covering numbers we are going to use in this section.\  
This section is divided in two subsections, in the Subsection \ref{estimates1} we present the detailed proof of Theorem \ref{coveringestimatesgeo}, and since the proof of Theorem \ref{coveringestimateshar} is very similar, we present very briefly its proof in Subsection \ref{estimates2}.

The covering numbers as defined in Section \ref{intro} is the particular case of the following.\ Consider $(X, \|\cdot\|_X), (Y,\|\cdot\|_Y)$ Banach spaces and $B_X$ the unit ball in $X$, then the \emph{covering numbers} of an operator $T: X \longrightarrow Y$ are given by
$$
\mathcal{C} (\epsilon, T):= \mathcal{C}(\epsilon,T(B_X)), \quad \epsilon>0.
$$

It is not hard to see that if $\| T\| \leq \epsilon $, then $\mathcal{C}(\epsilon, T)=1$, for any $\epsilon>0$, and that $C(\epsilon, T)$ is non increasing in $\epsilon$.\ The following properties can be found in \cite{Carl,pietsch}, for example.

If $S,\, T: X\longrightarrow Y$ and $R: Z \longrightarrow X$ are operators on real Banach spaces, then
\begin{enumerate}
\item[c1.] $\mathcal {C} (\epsilon + \delta, T+S) \leq \mathcal{C} (\epsilon, T) \hspace{0.1cm} \mathcal{C}(\delta, S)$, and
 \item[c2.] $\mathcal{C}(\epsilon \delta,TR)\leq \mathcal{C}(\epsilon, T)\,\mathcal{C}(\delta, R)$, for any $\epsilon, \delta>0$;
\item[c3.] if $n:=\mbox{rank}(T) <\infty$, then $\mathcal{C}(\epsilon, T)\leq \left( 1+2\| T\|/\epsilon \right)^ n$, for any $\epsilon>0$.
\end{enumerate}
For the lower estimates we will need the following inequality for the covering numbers for two $n$-dimensional Hilbert spaces $X$ and $Y$, 
\begin{equation}\label{CN1}
\left|\det \sqrt{T^{*}T}\right| \left(\frac{1}{\epsilon}\right) ^n \leq \mathcal{C}(\epsilon,T), \quad \epsilon>0.
\end{equation}

\subsection{Proof of Theorem \ref{coveringestimatesgeo}}\label{estimates1}

The proof will be an immediate consequence of the upper estimate obtained in Proposition \ref{thmuppergeo} and the lower estimate obtained in Proposition \ref{thmlowergeo} bellow.\ 

In order to present the proof we fix the following short notation for some quantities obtained in Proposition \ref{thmNorms},
\begin{equation}\label{kappas}
\kappa= \|I_K\|, \quad  \kappa_m= \|P_m\|, \quad \mbox{and} \quad \kappa_m^s= \|P_m^s \|, \quad \mbox{for $m=1,2,\ldots$}.
\end{equation}

We observe that if $K$ an isotropic positive definite kernel represented by the series expansion (\ref{WKlegendre}), with $a_0\neq 0$, and there exists $0<\theta<1$ such that  $a_k \leq \theta\, a_{k-1}$, for $k=1,2\ldots$, then Proposition \ref{thmNorms} implies 
\begin{equation}\label{CS}
0<(k_m^s)^2 \leq \frac{\theta^{m+1}}{1-\theta} \, a_0, \quad m=1,2, \ldots.
\end{equation} 

\begin{prop}\label{thmuppergeo}  Let $K$ an isotropic positive definite kernel represented by the series expansion (\ref{WKlegendre}).\ If there exists $0<\theta<1$ such that  $a_k \leq \theta\, a_{k-1}$, for $k=1,2\ldots$, then 
$$
\limsup_{\epsilon\to 0^+}\frac{\ln (\mathcal{C} (\epsilon, I_K))}{[\ln(1/\epsilon)]^{d+1}}\leq \frac{2^{d+2}}{d! \left({\ln\left( {1}/{\theta}\right)}\right)^d}. 
$$
\end{prop}

\pf 
With no loss of generality we can consider $a_0\neq 0$.\ Otherwise, in the proof bellow, we just replace $a_0$ by the first positive term of the sequence $\{a_n\}$.

For $\epsilon >0$, we consider $m:=m(\epsilon)\in \mathbb{Z}_+$ as follows
\begin{equation}\label{alphaMEpsilonexp}
     \left(\frac{\theta^{m+1}}{1-\theta} \,a_0\right)^{1/2}    < \frac{\epsilon}{2} < \left(\frac{\theta^{m}}{1-\theta}\, a_0\right)^{1/2}.
\end{equation} 
Applying the logarithm above, we have that $m< \ln (\epsilon^2 (1-\theta) / 4\, a_0 )/\ln\theta$, and since 
$$\frac{\ln (\epsilon^2 (1-\theta) / 4\, a_0 )}{\ln\theta} \approx 2\frac{\ln(1/\epsilon)}{\ln(1/\theta)},
$$
for $\epsilon>0$ sufficiently small, or equivalently, for $m\in\mathbb{Z}_+$ sufficiently large, the following estimate holds
\begin{equation}\label{mEpsilonexp}
m< 2\frac{\ln(1/\epsilon)}{\ln(1/\theta)}, \quad m=m_0,m_0+1,\ldots.
\end{equation}
It is important to notice that considerations above, inequality (\ref{CS}), and inequality (\ref{alphaMEpsilonexp}) imply 
$$
\mathcal{C}\left(\epsilon , P_m^s \right)=1, \quad m=1,2,\ldots. 
$$
Thus, properties c1. and c3. for covering numbers, respectively, lead us to 
$$
\mathcal{C} (\epsilon, I_K = P_m+P_m^s) \leq \left(1+ 4\| P_m \|/\epsilon \right)^ {\mbox{rank}(P_m)} \leq \left(1+ 4\kappa_m/\epsilon \right)^ {\dim{\mathcal{V}_m}},
$$
and an application of Corollary \ref{dimVm}, implies that the following inequality holds
$$
\mathcal{C} (\epsilon, I_K) < \left(1+ 4\kappa_m/\epsilon \right)^{4m^d/d!},
$$
for $m$ large, and we assume that it holds for $m\geq m_0$.\ Taking in account that $\kappa_m\leq\kappa$, for any $m$, and applying the logarithm in both sides of inequality above, we obtain
\begin{eqnarray*}
    \ln (\mathcal{C}(\epsilon, I_K ))< \frac{4m^d}{d!} \ln \left(1+4\kappa/\epsilon\right) \leq \frac{4m^d}{d!} \ln \left(8\kappa/\epsilon\right), \quad m=m_0,m_0+1,\ldots.
    \end{eqnarray*}    
Finally, inequality (\ref{mEpsilonexp}) implies
$$
\ln (\mathcal{C}(\epsilon, I_K ))< \frac{4m^d}{d!} \ln \left(8\kappa/\epsilon\right) \leq \frac{4}{d!}\left(\frac{2\ln (1/\epsilon)}{\ln\left(1/ \theta\right)}\right)^d  \ln \left( \frac{1}{\epsilon}\right),
$$ 
for $m=m(\epsilon)=m_0,m_0+1,\ldots$, and the proof follows. \eop

\begin{prop}\label{thmlowergeo} Let $K$ an isotropic positive definite kernel represented by the series expansion (\ref{WKlegendre}).\ If there exists $0<\delta<1$ such that  $a_k \geq \delta\, a_{k-1}>0$, for $k=1,2\ldots$, then 
\[
\frac{ (2d)^d }{d! (d+\ln(1/\delta)-1)^d (d+1)^{d+1}}\leq \liminf_{\epsilon\to 0^+}\frac{\ln(\mathcal{C}(\epsilon, I_K ))}{[\ln(1/\epsilon)]^{d+1}} .
\]
\end{prop}

\pf For $m=1,2,\ldots$ we consider the following composition operator
$$
T_m : \mathcal{V}_m \stackrel{ J_m }{\longrightarrow} \mathcal{H}_K  \stackrel{I_K}{\longrightarrow} C(\mathbb{S}^d) \stackrel{{L_m}}{\longrightarrow}  L^2 ( \mathbb{S}^d ) \stackrel{P_m }{\longrightarrow} L_m\circ I_K\circ J_m(\mathcal{V}_m),
$$
where $J_m$ stands for the identity operator given by $\mathcal{V}_m\hookrightarrow\mathcal{H}_K$, $L_m$ and
$P_m $ are the orthogonal projections on $L_m I_K J_m(\mathcal{V}_m)$.\ This operator induces $T^{*}_m\,T_m$ and the representing matrix satisfies
\begin{equation}\label{det}
\det(T_m ^* T_m)= \prod_{s=0}^{m} \left(\frac{a_s}{\tau_{d}^{s}}\right)^{\tau_{d}^{s}}.
\end{equation}
It is clear that $\|P_m\|=\|J_m\|=\Vert L_m \Vert=1$.\ 

For any $\epsilon>0$, properties c1. and c2, for the covering numbers, imply 
\begin{eqnarray*}
\mathcal{C} (\,\epsilon, T_m ) & = &\mathcal{C} (\,\epsilon, P_m L_m I_K J_m ) \\ & \leq & \mathcal{C} (1, P_m ) \mathcal{C} (1,L_m ) \mathcal{C} (\epsilon, I_K )\mathcal{C} (1, J_ m )\\ & = & \mathcal{C} (\epsilon, I_K ), \quad m=1,2,\ldots.
\end{eqnarray*}
The lower bound for the covering number in inequality (\ref{CN1}), leads us to
$$
\sqrt{ \det (T_m ^* T_m )}\left(\frac{1}{{}\,\epsilon}\right)^{\dim \mathcal{V}_m} \leq \mathcal{C} (\,\epsilon, T_m  ) \leq \mathcal{C} (\epsilon, I_K ), 
$$
and due Corollary \ref{dimVm} and formula (\ref{det}), we obtain
$$
\sqrt{ \det (T_m ^* T_m )}\left(\frac{1}{{}\,\epsilon}\right)^{\dim \mathcal{V}_m} = \prod_{k=0}^{m} \left(\frac{a_k}{\tau_d^k}\right)^{\tau_d^k/2}\left(\frac{1}{\,\epsilon}\right)^{\tau_{d+1}^{m}}, \quad \epsilon>0, \quad m=1,2,\ldots.
$$
Thus,
$$
\prod_{k=0}^{m} \left(\frac{a_k}{\tau_d^k}\right)^{\tau_d^{k}/2}\left(\frac{1}{\,\epsilon}\right)^{\tau_{d+1}^{m}}\leq \mathcal{C} (\epsilon, I_K ), \quad \epsilon>0, \quad m=1,2,\ldots. 
$$
Applying the logarithm in both sides of inequality above, we conclude that for any $\epsilon>0$, 
\begin{equation}\label{lower1geo}
     \frac{1}{2}\sum_{k=0}^{m} \tau_d^{k}\ln \left(a_k/\tau_d^{k}\right) - \tau_{d+1}^{m} \ln({}\,\epsilon)\leq \ln(\mathcal{C}(\epsilon, I_K )), \quad m=1,2,\ldots.
\end{equation}
We observe that, for $m=1,2,\ldots$, the estimate 
$$
\delta^{m}a_0/\tau_{d}^{m}\leq a_m/\tau_{d}^{m} \leq a_k/\tau_d^{k}, \quad k=0,1,\ldots m,
$$
implies 
\begin{equation}\label{lower2geo}
\frac{\tau_{d+1}^{m} }{2} \left( \ln (\delta^{m}a_0)-\ln \tau_{d}^{m} \right) \leq \frac{1}{2}\sum_{k=0}^{m} \tau_d^{k} \ln\left(a_k/\tau_d^{k} \right), \quad m=1,2,\ldots.
\end{equation}
If we write
$$
I_m:=\frac{\tau_{d+1}^{m} }{2} \left( \ln (\delta^{m}a_0)-\ln \tau_{d}^{m} \right), \quad m=1,2,\ldots,
$$
then Corollary \ref{dimVm} ensures that 
\begin{eqnarray*}
 \frac{m^d}{2 d!}\left[\ln (\delta^{m}a_0) -m(d-1) - \ln (4/(d-1)!)\right]< I_m, \quad m=m_0, m_0+1,\ldots, 
\end{eqnarray*}
for some $m_0\in\mathbb{Z}_+$.\ This estimate and inequalities (\ref{lower1geo}) and (\ref{lower2geo}) imply
$$
\frac{m^d}{2d!}\left[m\ln \delta + \ln a_0 -m(d-1) - \ln (4/(d-1)!) - 2\ln(\epsilon)\right]< \ln(\mathcal{C}(\epsilon, I_K )),
$$
for any $\epsilon>0$ and $m=m_0,m_0+1,\ldots$.

The final steps of the proof is to solve the optimization problem given by the following quantity,
$$ 
\varphi_\epsilon (m):= -m^d [ \ln (4 \,\epsilon^2 /a_0  (d-1)!)]-m^{d+1}[\ln(1/\delta)+(d-1)],
$$
with $\epsilon>0$ and $m=m(\epsilon)$ such that $1/(m+1)\leq\epsilon <1/m$. \ Then, we obtain implies
\begin{equation}\label{lower5geo}
\frac{\varphi_\epsilon (m)}{ 2d!}< \ln(\mathcal{C}(\epsilon, I_K )), \quad m=m_0, m_0+1, \ldots,
\end{equation}
The critical point of $\varphi_{\epsilon}$, as a function on $\mathbb{R}$, is given by
$$
c=-\frac{d}{(d+1)[d+\ln(1/\delta)-1]} \ln(4 \,\epsilon^2 /a_0 (d-1)!).
$$
Since $\lfloor c \rfloor \approx  c \approx \lceil c \rceil $, we can estimate (\ref{lower5geo}) in terms of $\varphi_{\epsilon}(c)$, even if $c$ is not an integer (and it my not be).\ A simple calculation leads us to 
\begin{eqnarray*}
\frac{\varphi_\epsilon (c)}{2{}d!} =\frac{[- \ln(4 \,\epsilon^2 /a_0 {}(d-1)!)]^{d+1}}{2{}d! \,[d+\ln(1/\delta)-1]^d} \frac{d^d}{(d+1)^{d+1}},
 \end{eqnarray*}
we can consider $a_0$ such that $a_0(d-1)!/4\geq 1$, and then inequality (\ref{lower5geo}) implies
$$
\frac{ (2d)^d }{{}d! (d+\ln(1/\delta)-1)^d (d+1)^{d+1}} < \frac{\ln(\mathcal{C}(\epsilon, I_K ))}{(\ln(1/\epsilon))^{d+1}},
$$
for $1/(m+1)\leq\epsilon <1/m$ and $m=m_0, m_0+1, \ldots$. The proof follows.
\eop

\subsection{Proof of Theorem \ref{coveringestimateshar}}\label{estimates2}

The proof of Theorem \ref{coveringestimateshar} will be an immediate consequence of the upper estimates (Proposition \ref{thmupperhar}) and the lower estimates (Proposition \ref{thmlowerhar}).\ The proofs are similar to the previous case and we will short the presentation for the similar details.

We observe that if $K$ is an isotropic positive definite kernel represented by the series (\ref{WKlegendre}) with sequence of coefficients satisfying
$$\limsup_{n\to\infty}a_n/n^{-\beta}\leq c,$$ for some $\beta>1$ and $c>0$, then there exists $m_0\in\mathbb{N}$ such that
\begin{equation}\label{CSpol}
0< (k_m^s)^2 \leq \frac{c}{\beta-1}{(m+1)}^{-\beta+1}, \quad m=m_0,m_0+1, \ldots,
\end{equation} 
with $\kappa_m^s$ as introduced in (\ref{kappas}).\ 

\begin{prop}\label{thmupperhar} Let $K$ an isotropic positive definite kernel given by the series expansion (\ref{WKlegendre}).\ If there exists $\beta>0$ and $c_1>0$ such that $a_k/ k^{-d-\beta} \leq c_1$, for $k=n_0, n_0+1, \ldots$, then 
$$
\limsup_{\epsilon\to 0^+}\frac{\ln (\mathcal{C} (\epsilon, I_K))}{(1/\epsilon)^{2d/(\beta+d-1)}\ln (1/\epsilon)}\leq \frac{4}{d!} \left(\frac{4 c_1 }{\beta+d-1}\right)^{d/(\beta+d-1)}. 
$$
\end{prop}
\pf For $\epsilon>0$, we consider $m:=m(\epsilon)\in \mathbb{N}$ such that
\begin{equation}\label{alphaMEpsilonpol}
     \left[\frac{c_1}{\beta+d-1}(m+1)^{-\beta-d+1}\right]^{1/2}    < \frac{\epsilon}{2 }< \left[\frac{c_1}{\beta+d-1}m^{-\beta-d+1}\right]^{1/2}.
\end{equation} 
In this case, we have
\begin{equation}\label{mEpsilonpol}
m(\epsilon) \approx \left(\frac{4 c_1 }{\beta+d-1}\right)^{1/(\beta+d-1)} \left( \frac{1}{\epsilon}\right)^{2/(\beta+d-1)}, \quad \mbox{as $\epsilon\to 0^+.$}
\end{equation}
Repeating exactly the same steps in the proof of Proposition \ref{thmuppergeo}, we obtain
\begin{eqnarray*}
    \ln (\mathcal{C}(\epsilon, I_K ))< \frac{4m^d}{d!} \ln \left(8\kappa/\epsilon\right), \quad m=m_0, m_0+1, \ldots,
    \end{eqnarray*}
for some $m_0\geq n_0$. An application of the strong equivalence (\ref{mEpsilonpol}) now, implies
$$
m^d\ln \left(8\kappa/\epsilon\right)
    < \left(\frac{4\,c_1 }{\beta+d-1}\right)^{d/(\beta+d-1)} \left( \frac{1}{\epsilon}\right)^{2d/(\beta+d-1)} \ln \left( \frac{1}{\epsilon}\right), \quad q_0, q_0+1, \ldots
$$
for some $q_0\geq m_0$. The proof follows. \eop

\begin{prop}\label{thmlowerhar} Let $K$ an isotropic positive definite kernel given by the series expansion (\ref{WKlegendre}) with $\{a_n\}$ a non-increasing sequence.\ If there exist $\rho>1$ and $c_2>0$ such that $a_k/k^{-\rho} \geq c_2$, for $k= 1,2,\ldots$, then
$$ 
\ln{\sqrt{a_0}}+ \frac{2(\rho+d-1)}{e \, d\, d!}  \left(\frac{c_2}{2 } \right)^{d/(\rho+d-1)} \leq \liminf_{\epsilon\to 0}\frac{\ln(\mathcal{C}(\epsilon, I_K ))}{(1/\epsilon)^{d/2(\rho+d-1)}} 
$$
\end{prop}

\pf Following the same steps of the proof of Proposition \ref{thmlowergeo}, from formula (\ref{det}) to formula (\ref{lower1geo}), for any $\epsilon>0$, we attain the following lower bound

\begin{equation}\label{lower1}
J_m:=  \frac{1}{2}\sum_{k=0}^{m} \tau_d^{k}\ln \left(a_k/\tau_d^{k}\right) - \tau_{d+1}^{m} \ln({\sqrt{\epsilon}})\leq \ln(\mathcal{C}(\epsilon, I_K )), \quad m=1,2,\ldots.
\end{equation}

For $m=1,2,\ldots$, we write 
\begin{equation}
I_m:=\frac{1}{2}\left(\ln a_0 + \sum_{k=1}^{m}\tau_d^k\ln a_k + \frac{1}{2}\sum_{k=1}^{m}\tau_d^k\ln(1/\tau_d^k)\right)= \frac{1}{2}\sum_{k=0}^{m} \tau_d^{k}\ln \left(a_k/\tau_d^{k}\right),
\end{equation}
and if $a_0\geq 1$ we observe that, for $m=1,2,\ldots$, 
\begin{eqnarray*}
 I_m &\geq & \frac{1}{2}\left(\ln a_0 + (\tau_{d+1}^m-1)\ln c_2 -\rho(\tau_{d+1}^m-1)\ln m+\sum_{k=1}^{m}\tau_d^k\ln(1/\tau_d^k)\right).
\end{eqnarray*}
For $m=1,2, \ldots$ and $k=1,2\ldots,m$, formula (\ref{dim}) implies $\tau_d^k\leq 2k^{d-1} \leq 2m^{d-1}$, then we gain
$$
\ln(1/\tau_d^k)\geq \ln (1/2m^{d-1}),
$$
and, consequently 
\begin{eqnarray*}
 I_m &\geq & \frac{1}{2}\left[\ln a_0 + (\tau_{d+1}^m-1)\ln c_2 -\rho(\tau_{d+1}^m-1)\ln m+ (\tau_{d+1}^m-1)\ln (1/2m^{d-1})\right]
 \\ &=& \ln \sqrt{a_0}+ \frac{(\tau_{d+1}^m-1)}{2}\left[\ln c_2 -\rho\ln m+\ln (1/2m^{d-1})\right]\\ &=& \ln \sqrt{a_0}+ \frac{(\tau_{d+1}^m-1)}{2}\left[\ln c_2 -(\rho-1)\ln m+\ln (1/2m^d)\right].
\end{eqnarray*}
From inequality (\ref{lower1}), we have 
\begin{equation}\label{Jm}
J_m=I_m- \tau_{d+1}^m\ln(\epsilon)\leq\ln(\mathcal{C}(\epsilon, I_K )), \quad m=1,2,\ldots.
\end{equation}
An application of (\ref{dimapproxm}) implies 
$$
\frac{\tau_{d+1}^m-1}{2}<\frac{{2} m^d}{d!} \quad \mbox{and} \quad -\tau_{d+1}^m<-\frac{m^d}{d!}, \quad m=m_0, m_0+1, \ldots,
$$
then 
$$
\ln \sqrt{a_0}+ \frac{{2} m^d}{d!}\left[\ln c_2 -(\rho-1)\ln m+\ln (1/2m^d)\right] -\frac{m^d}{d!}\ln(\epsilon)<J_m.
$$
Therefore,
\begin{equation}\label{base}
\ln \sqrt{a_0}+ \frac{{2}m^d}{d!}\left[\ln (c_2/\epsilon^{1/{2}}) - (\rho-1)\ln m -\ln(2m^{d})\right]<J_m
\end{equation}
for $m=m_0, m_0+1, \ldots$.

In order to proceed analogously to the final steps of the proof of the Proposition \ref{thmlowergeo}, we consider $\epsilon>0$ and $m=m(\epsilon)$ such that $1/(m+1)\leq \epsilon <1/m$, and
$$ 
\varphi_{\epsilon} (m):=  m^d\left[\ln (c_2/2\epsilon^{1/{2}}) - (\rho-1)\ln m -\ln(m^d)\right].
$$
Due inequality (\ref{Jm}), it holds 
\begin{equation}\label{lower3}
\ln{\sqrt{a_0}}+\frac{2\varphi_\epsilon (m)}{d!}< \ln(\mathcal{C}(\epsilon, I_K )), \quad m=m_0, m_0+1, \ldots,
\end{equation}
We observe that from inequality (\ref{lower3}), we will have 
\begin{equation}\label{lower5}
\ln{\sqrt{a_0}}+\frac{2 \varphi_{\epsilon} (c)}{  d!}< \ln(\mathcal{C}(\epsilon, I_K )),
\end{equation}
for the critical point $c$ of $\varphi_{\epsilon}$.\ In this case the critical point is 
$$
c= \left(\frac{c_2}{2 \,\epsilon^{1/2} }\right)^{1/(\rho+d-1)} e^{-1/d}, 
$$ and we have 
\begin{eqnarray*}
   \ln{\sqrt{a_0}}+ \frac{2\varphi_\epsilon (c)}{ d!} = \ln{\sqrt{a_0}}+ \frac{2(\rho+d-1)}{e\, d\, d!}  \left(\frac{c_2}{2 } \right)^{d/(\rho+d-1)} \left(\frac{1}{\epsilon}\right)^{d/2(\rho+d-1)},
\end{eqnarray*}
and the first estimate follows. \eop

\section{Asymptotic estimates for the covering numbers}\label{consequences}

In this section we consider the strong equivalence stated in formula (\ref{dimapprox}) in order to explore the problem of finding the asymptotic constant for the covering numbers of kernels as in Theorem \ref{coveringestimatesgeo}.\ This subject is explored in \cite{kuhn22} in a more abstract setting.\  

\begin{prop}\label{propasympconstant}  Let $K$ an isotropic positive definite kernel represented by the series expansion (\ref{WKlegendre}).\ The following hold
\begin{itemize}
\item[I)] If there exists $0<\theta<1$ such that  $a_k \leq \theta\, a_{k-1}$, for $k=1,2\ldots$, then 
$$
\limsup_{\epsilon\to 0^+}\frac{\ln (\mathcal{C} (\epsilon, I_K))}{[\ln(2a_0^{1/2}/\epsilon(1-\theta)^{1/2})]^{d+1}}\leq \frac{2^{d+1}}{{}d! \left({\ln\left( {1}/{\theta}\right)}\right)^d}. 
$$
\item[II)] If there exists $0<\delta<1$ such that  $a_k \geq \delta\, a_{k-1}>0$, for $k=1,2\ldots$, then 
\[
\frac{ 2^{d+1} d^d  }{{}d! (d+\ln(1/\delta)-1)^d (d+1)^{d+1}}\leq \liminf_{\epsilon\to 0^+}\frac{\ln(\mathcal{C}(\epsilon, I_K ))}{[\ln(\sqrt{(a_0(d-1)!)/2})/\epsilon)]^{d+1}} .
\]
\end{itemize}
\end{prop}

\pf By proceeding analogously as in the proof of Proposition \ref{thmuppergeo}, for $\epsilon >0$, we consider $m:=m(\epsilon)\in \mathbb{Z}_+$ as in inequality (\ref{alphaMEpsilonexp}).\ In this case, we notice that 
\begin{equation}\label{m1}
m\approx \frac{\ln (\epsilon^2 (1-\theta) / 4\, a_0 )}{\ln\theta}= \frac{2\ln(2a_0^{1/2}/\epsilon(1-\theta)^{1/2})}{\ln(1/\theta)},    
\end{equation}
and, by an application of equivalence (\ref{dimapprox}), we have
$$
L_m:=\left(1+ 4\kappa_m/\epsilon \right)^{\dim\mathcal{V}_m}\approx\left(1+ 4\kappa/\epsilon \right)^{2m^d/{}d!}.
$$
It follows that
$$
\ln L_m\approx \frac{2m^d}{{}d!} \ln \left(1+4\kappa/\epsilon\right),
$$
and we have that $\mathcal{C} (\epsilon, I_K) \leq \left(1+ 4\kappa_m/\epsilon \right)^{\dim\mathcal{V}_m}=L_m$ holds true, then equivalence (\ref{m1}) implies

\begin{eqnarray*}
\limsup_{\epsilon\to 0^+}\frac{\ln (\mathcal{C}(\epsilon, I_K ))}{[\ln(2a_0^{1/2}/\epsilon(1-\theta)^{1/2})]^{d+1}} & \leq & \frac{2^{d+1}}{{}d!(\ln\left(1/ \theta\right))^{d}},   
\end{eqnarray*}
and the proof of I) follows.

For the proof of II), we proceed analogously as in the proof of Proposition \ref{thmlowergeo}.\ For $\epsilon >0$, we have
\begin{equation}\label{m0}
I_m:= \frac{1}{2}\sum_{k=0}^{m} \tau_d^{k}\ln \left(a_k/\tau_d^{k}\right) - \tau_{d+1}^{m} \ln(\epsilon)\leq \ln(\mathcal{C}(\epsilon, I_K )), \quad m=1,2,\ldots.
\end{equation}
If we write
$$
J_m:=\frac{m^d}{{}d!}\left[m\ln \delta + \ln a_0 -m(d-1) - \ln (2/{}(d-1)!) - 2\ln(\epsilon)\right],
$$
then $I_m\approx J_m$, by an application of equivalence (\ref{dimapprox}).\ We consider
$$
\varphi_\epsilon (m):= -m^d [ \ln (2 \,\epsilon^2 /a_0 \,{} (d-1)!)]-m^{d+1}[\ln(1/\delta)+(d-1)],
$$
with $\epsilon>0$ and $m=m(\epsilon)$ such that $1/(m+1)\leq\epsilon <1/m$.\ Thus,
\begin{eqnarray*}
\frac{\varphi_\epsilon (c)}{{}d!} =\frac{[- \ln(2 \,\epsilon^2 /a_0 {}(d-1)!)]^{d+1}}{{}d! \,[d+\ln(1/\delta)-1]^d} \frac{d^d}{(d+1)^{d+1}}\leq J_m,
 \end{eqnarray*}
with $c$ the critical point of $\varphi_{\epsilon}$.\ Since,
$$
\ln(2 \,\epsilon^2 /a_0 {}(d-1)!)= -2\ln(\sqrt{(a_0(d-1)!)/2})/\epsilon),
$$ 
inequality above, and inequality (\ref{m0}) imply the inequality stated in II). \eop

\begin{cor} \label{sharper} Let $K$ an isotropic positive definite kernel represented by the series expansion (\ref{WKlegendre}).\ If  $a_0=1$ and $a_k \leq \theta^k$, for some $0<\theta<1$ and $k=1,2,\ldots$, then 
$$
\limsup_{\epsilon\to 0^+}\frac{\ln (\mathcal{C} (\epsilon, I_K))}{[\ln(1/\epsilon)]^{d+1}}\leq \frac{2^{d+1}[\ln(2e/(1-\theta)^{1/2})]^{d+1}}{{}d! \left({\ln\left( {1}/{\theta}\right)}\right)^d}. 
$$
\end{cor}
 
\pf The inequality stated easily follows from
\begin{eqnarray*}
\limsup_{\epsilon\to 0^+}\frac{\ln (\mathcal{C} (\epsilon, I_K))}{[\ln(1/\epsilon)]^{d+1}}& \leq & [\ln(2e/(1-\theta)^{1/2})]^{d+1}\limsup_{\epsilon\to 0^+}\frac{\ln (\mathcal{C} (\epsilon, I_K))}{[\ln(2/\epsilon(1-\theta)^{1/2})]^{d+1}},
\end{eqnarray*}
and Proposition \ref{propasympconstant}.\ The inequality above is due the fact that
\begin{eqnarray*}
\frac{\ln (1/\epsilon)}{\ln(2/\epsilon(1-\theta)^{1/2})}>\frac{1}{\ln(2e/(1-\theta)^{1/2})}, 
\end{eqnarray*}
for $0<\epsilon< e$. \eop

The following is a direct consequence of item II) in Proposition \ref{propasympconstant} and Corollary \ref{sharper}, since the upper and the lower bounds became equal under the assumptions.   

\begin{cor}\label{thmasympconstant}  Let $K$ an isotropic positive definite kernel represented by the series expansion (\ref{WKlegendre}) with the sequence of coefficients satisfying $a_0=1$ and $a_k \leq \theta a_{k-1}$, for $k=1,2,\ldots$ for some $0<\theta<1$.\ If $\{a_k\}_k$ satisfies $a_k\geq \delta_{\theta, d} a_{k-1}$ for
$$
\delta_{\theta, d}:=\theta^{c(\theta,d)} , \quad \mbox{with}\quad c(\theta,d):=d/\left[\ln\left(\frac{2e}{(1-\theta)^{1/2}}\right)^{d+1}\right]^{1+1/d}+(d-1)/\ln \theta, 
$$
then
$$
\lim_{\epsilon\to 0^+}\frac{\ln (\mathcal{C} (\epsilon, I_K))}{[\ln(1/\epsilon)]^{d+1}}=\frac{2^{d+1}[\ln(2e/(1-\theta)^{1/2})]^{d+1}}{{}d! \left({\ln\left( {1}/{\theta}\right)}\right)^d}. 
$$
\end{cor}

A consequence of Theorem \ref{coveringestimatesgeo} is the following well known asymptotic behavior for the covering numbers of the Gaussian like kernel (see \cite{jordao,Kuhn} and references therein).

\begin{cor} Let $K$ be an isotropic positive definite kernel given by the series expansion (\ref{WKlegendre}).\ If there exists $0<\theta<1$ such that  $a_k\asymp \theta^k$, then
\[\ln(\mathcal{C}(\epsilon, I_W ))\asymp \left[ \ln \left(1/\epsilon \right) \right]^{d+1},  \quad \mbox{as}\quad \epsilon \longrightarrow 0.\]
\end{cor}

As a consequence of Theorem \ref{coveringestimateshar}, we obtain an improvement of Theorem 3.1 in \cite{Azevedo}, for the particular case of the $d$-dimensional sphere as the compact two-point homogeneous space.\ The results in this paper can be easily extended to the the framework considered in \cite{Azevedo} making the necessary adaptations on the quantities presented in Section \ref{RKHS}.

\begin{cor} Let $K$ an isotropic positive definite kernel given by the series expansion (\ref{WKlegendre}).\ If there exists $\gamma>1$ such that $a_k \asymp k^{-d-\gamma}$, then there exist positive constants $A, B$ do not depending on $\epsilon$, such that
$$
A\left(\frac{1}{\epsilon}\right)^{d/2(d+\gamma-1)}\leq \ln(\mathcal{C}(\epsilon, I_K ))\leq B \left( \frac{1}{\epsilon}\right)^{2d/(\gamma+d-1)} \left[\ln \left( \frac{1}{\epsilon}\right)\right].
$$
\end{cor}

\pf The proof for the lower bound is an direct application of Proposition \ref{thmlowerhar} for $\rho = \gamma$, and the upper bound from Proposition \ref{thmuppergeo} for $\beta = d+\gamma$.\eop

\section{The basic decay of the Fourier coefficients of isotropic positive definite kernels}

The kernel $K$ as defined in (\ref{WKlegendre}), can be easily fitted in the context considered in articles \cite{carrijo,jordao,jordao1}. If $d\geq 5$ is odd, then an immediate consequence of \cite[Theorem 3.1]{jordao1} or in \cite[Theorem 4]{carrijo}, for the sequence $\{a_k\}$ is the following decay $a_k=O(k^{-d})$. Therefore, it is very natural the assumption $a_k/k^{d+\beta}=O(1)$, for $\beta>0$, considered in the Theorem 1.2.

In fact, we consider the shifting operator defined by
$$
S_t f(x)=\frac{1}{R_d(t)} \int_{R_x^t} f(y) d \sigma_r(y), \quad x \in S^d, \quad f \in L^1\left(S^d\right), \quad t \in(0, \pi),
$$
in which $d \sigma_r(y)$ is the volume element of the rim $R_x^t:=\left\{y \in S^d: x \cdot y=\cos t\right\}$ and $R_d(t)=\sigma_{d-1}(\sin t)^{d-1}$ is the total volume. It is well known that $\|S_t\| \leq 1$ and $\{S_t\}$ is a family of multipliers operators with $\{P_k^{(m-1) / 2}(\cos t)\}$ the sequence of multipliers, i.e.,
$$
S_t(f) \sim \sum_{k=0}^{\infty} P_k^{(m-1) / 2}(\cos t)P_k(f), \quad f\in L^2\left(S^d\right), \quad t\in (0,\pi),
$$
with $P_k(f) = \sum_{j=1}^{\tau_d^k}\left\langle f, S_{k, j}\right\rangle_2 S_{k, j}$ as in Proposition \ref{thmNorms}, for each $k=0,1,2,\ldots$.\ For for $r>0$ an integes number, the generalized shifting operator (\cite{carrijo}) is
$$
S_{r, t}(f)=\frac{-2}{{2r \choose r}}  \sum_{j=1}^r(-1)^j{2r\choose r-j} S_{jt}(f), \quad t>0,
$$
it means that $S_{r,t}(f) \sim \sum_{k=0}^{\infty}m_r(k, t)P_k(f)$, for $f\in L^2(S^d)$, and $t\in (0,\pi)$, and the multiplier sequence is
$$
m_r(k, t)=\frac{-2}{{2r \choose r}}\sum_{j=1}^r(-1)^j{2r\choose r-j} P_k^d(\cos (j t)), \quad k=0,1 \ldots, \quad t>0,
$$
satisfying the asymptotic relation $\left(m_r(k, t)-1\right)\asymp [\min\{1,kt\}]^{2r}$ (\cite[Lemma 11]{carrijo}). Thus, it is not hard to derive the following
\[
I_t(x,y):=S_{r,t}(K(x,\,\cdot\,))(y)- K(x, y) = \sum_{k=0}^\infty \frac{a_k}{\tau_d^k} \left(m_r(k, t)-1\right)\sum_{j=1}^{\tau_d^k}S_{k,j}(x)S_{k,j}(y), \quad x,y\in S^d.
\]
An application of the Cauchy-Schwarz inequality for square summable sequences implies
\[
|I_t(x,y)| \leq \left(\sum_{k=0}^\infty \sum_{j=1}^{\tau_d^k}\frac{a_k}{\tau_d^k}\left(m_r(k, t)-1\right)^2|S_{k,j}(x)|^2\right)^{1/2} \left(\sum_{k=0}^\infty \sum_{j=1}^{\tau_d^k}\frac{a_k}{\tau_d^k}|S_{k,j}(y)|^2\right)^{1/2},
\]
for any $t\in (0,\pi)$, and $x,y\in S^d$.\ The standard estimate $|S_{k,j}(y)|^2\leq \tau_d^k$ (\cite[Theorem 2.29]{morimoto}) and $\left(m_r(k, t)-1\right)\asymp [\min\{1,kt\}]^{2r}$, lead us to 
\[
|I_t(x,y)| \leq c' t^{2r}\left(\sum_{k=n_0}^\infty a_k k^{4r}\right)^{1/2},
\]
for $c'= c\kappa$, with $c$ some positive constant and $\kappa$ as defined in (\ref{thmNorms}). 

The fact that isotropic positive definite kernels have continuous derivative of order $\lambda:=\lfloor(d-1)/2\rfloor$ (see \cite[Theorem 1.2]{ziegel}), then $\sum_{k=0}^{\infty}a_kk^{\lambda}<\infty$. Consequently,
$a_k=O(k^{-(\lambda+1)})$ and if $4r>\lambda+2$, we gain
$$
|I_t(x,y)|  \leq  c't^{2r}\left(\sum_{k=n_0}^\infty k^{4r-\lambda-1}\right)^{1/2} \leq t^{2r} \frac{c'}{(4r-\lambda -2)^{1/2}}.
$$
In this case, the kernel $K$ satisfies the H\"{o}lder condition defined in formula (2.5) in \cite{jordao} (or in formula (2.2) in \cite{jordao1}), for the family of multipliers operators $\{S_{r,t}\}$ and order $\rho =2r$, since we showed that
$$
|S_{r,t}(K(x,\,\cdot\,))(y)- K(x, y)|\leq C t^{\rho}, \quad x,y\in S^d,
$$
for some constant $C>0$. Thus, \cite[Corollary 3.2]{jordao1} or \cite[Corollary 5]{carrijo} promptly imply the basic decay $a_n/\tau_d^n=O(n^{-d-2r})$. The asymptotic relation (\ref{dimapprox}) implies $a_n=O(n^{-2r-1})$, for $4r>\lambda+2$. The previous condition is fulfilled for $r=\lambda$, since $d\geq 5$, and we gain $a_n=O(n^{-d})$.  

\vspace{0.25cm}

\noindent\textbf{Acknowledgements.} The second author was financed in part by Coordenação de Aperfeiçoamento de Pessoal de Nível Superior - Brasil (CAPES) - Finance Code 001, and she is grateful to the Professor Douglas Azevedo by the guidance to this interesting problem.\ The corresponding author is supported by the São Paulo Research Foundation (FAPESP - grant no. 2022/11032-0).\ 

The authors thank the anonymous referee for the careful reading of the article, for the comments, and the questions, guiding us to present this improved version of the paper.\ Specially, we thank the referee for pointing out the problem of finding the asymptotic constant for the covering numbers (subject explored in \cite{kuhn22}).\ The comments offered leaded us to present the results in Section 4, and the comments in Section 5 about the basic decay of the Fourier coefficients of isotropic positive definite kernels.


\medskip

\noindent \textsc{T. Jord\~{a}o}\\
Departament of Mathematics\\
Universidade de S\~{a}o Paulo\\
13566-590 S\~{a}o Paulo, Brazil\\
\textsc{Email}: tjordao@icmc.usp.br

\medskip

\noindent {\textsc{K. Gonzalez}}\\
Departament of Mathematics\\
Universidade de S\~{a}o Paulo\\
13566-590 S\~{a}o Paulo, Brazil\\
\textsc{Email}: karina.navarro@correounivalle.edu.co

\medskip

\noindent{\textsc{Keywords}:} isotropic kernel; covering numbers; Reproducing Kernel Hilbert Space

\smallskip

\noindent{\textsc{MSC}:} 47B06; 42C10; 46E22; 41A60

\end{document}